\documentclass[12pt]{article}
\usepackage{amsfonts}
\usepackage{amssymb}
\usepackage{amsmath}
\usepackage{amsthm}
\usepackage{graphicx}
\usepackage[latin5]{inputenc}
\begin{document}

\date{}

\author{\textbf{Celalettin Kaya, Hurşit Önsiper}}
\title{ON ORDINARITY AND LIFTING FIBRATIONS ON SURFACES}
\maketitle

\begin{abstract}

We investigate the relation between the ordinarity of a surface and of its
Picard scheme in connection with the problem of lifting fibrations of genus $g \ge 2$
on surfaces to characteristic zero.

\end{abstract}

We consider a separable, connected fibration $\pi : X \rightarrow
C$ with fiber genus $g \ge 2$ on a smooth projective surface $X$ 
over a field $k$ of characteristic $p$ and
we ask if $\pi$ lifts to characteristic zero. More precisely, we
look for a discrete valuation ring $R$ of mixed characteristic with
residue field $k$, a
smooth projective surface $\cal X$  and a smooth curve $\cal C$
over $R$ which give the following commutative diagram.

\begin{center}
$
\begin{array}{ccc}
  X & \hookrightarrow & {\cal X} \\
  \downarrow & & \downarrow \\
  C & \hookrightarrow & {\cal C} \\
  \downarrow & & \downarrow \\
  Spec(k) & \hookrightarrow & Spec(R) \\
\end{array}
$
\end{center}

Unless otherwise stated explicitly, we take $R = W(k)$. \\

Clearly, the problem makes sense only if $X$ lifts as a surface.
Then, since there is no obstruction to lifting the base curve $C$,
the problem becomes a special case of the general problem of
lifting morphisms to characteristic zero. \\

The effect of the Picard scheme $Pic_X$ on the lifting problem for $X$ appears
in the algebraization of the formal scheme $\hat X$, if exists. To verify that the surface $X$ lifts 
it suffices to prove that
an ample line bundle on $X$ lifts to $\hat X$; 
  the obstruction in this step is clearly  related to the structure of $Pic_X$.\\

There is yet a rather indirect relation between the properties of $Pic_X$ and the lifting problem for $X$.
This relation arises from the elementary observation that the failure of Bogomolov inequality is
an obstruction to lifting $X$. It was conjectured that Bogomolov inequality holds if
$Pic_X$ is smooth (\cite{P}). However, this conjecture was disproved by Jang (\cite{JJ}) by constructing
fibered surfaces with smooth generically ordinary fibers and with smooth Picard scheme $Pic_X$, for which
Bogomolov inequality does not hold. In particular these surfaces do not lift to characteristic zero. \\

For a smooth isotrivial fibration the lifting problem is naturally related to 
the $p$-rank of the fiber. This follows from a monodromy criterion for lifting (Theorem 1) 
combined with the relation between the structure of the automorphism group of a curve
and its $p$-rank (\cite{N}). As an immediate corollary we obtain the main result of the paper: \\

\textbf{Corollary 1 :}\textit{ Let $\pi : X \to C$ be a smooth ordinary fibration with fiber $F$ and monodromy group $G$.
Suppose that either $g(F) \le p-2$ or $G$ is cyclic. Then $\pi$ lifts to characteristic zero.} \\

This result is proved in the second section of the paper. In the first section of the paper we prove some basic results
concerning the ordinarity of a variety and of its Picard scheme. \\

We first fix our notation.

\begin{itemize}
    \item $k$ is a field (algebraically closed unless otherwise stated) of
characteristic $p > 0$.
    \item $R$ is a complete discrete valuation ring of characteristic zero with
residue field $k$.
    \item $W(k)$ is the ring of Witt vectors over $k$ and $W_n(k)=W(k)/m^n$ is the ring of Witt vectors of length n.
    \item $X$ is a projective smooth surface over $k$.
    \item A lifting of $X$ means a projective smooth scheme $\cal X$ over $S
= Spec(R)$ with special fiber ${\cal X} \times_S Spec(k) \cong X$.
\item $Pic_X$ is the Picard scheme of  $X$, $Pic^0_X$ is the connected component of $Pic_X$ containing identity.
\item $Alb_X$ is the Albanese variety of $X$.
\item $B_X^j = d(\Omega_X^{j-1})$.
\item $C$ is the Cartier operator.
\item $F, F_X$ denote the absolute and relative Frobenius morphisms, respectively.
\item $F$ also denotes the action of Frobenius on any cohomology group etc.
\item $\mu_p, \alpha_p$ are the standard infinitesimal group schemes.
\item For an abelian variety $A$,

$A[n]$ is the subgroup scheme \textit{kernel of multiplication by $n$},

$\tau(A) = ~\mbox{dim}_k (Hom(\alpha_p, A))$.

\end{itemize}

\textbf{I. Some basic results on ordinarity} \\

We recall the following definition. \\

\textbf{Definition 1}: A smooth proper variety $X$ is said to be \textit{ordinary } in dimension $m$ if
$H^i(X,B^j) = 0$ for all $i,j$ such that $i+j = m$. $X$ is \textit{ordinary}, if it is ordinary in all dimensions. \\

For an abelian variety $A$, Definition 1 coincides with the classical definition; $A$ is
	ordinary if $|A[p]| = p^n$ where $n = \mbox{dim}(A)$ (\cite{BlKato}, Example 7.4). In fact, one has the following well-known list
	of equivalent conditions which we include without proof for later use in the paper. \\

\textbf{Lemma 1 :} \textit{For an abelian variety $A$ of dimension $n$ the following are equivalent.}

 \textit{\begin{tabular}{l}
~a) $A$ is ordinary in the sense of Definition 1. \\
~b) $A$ is ordinary in dimension one. \\
~c) $A[p] = (\displaystyle {\mathbb Z}/{p \mathbb Z})^n \times \mu_p^n$. \\
~d) $\tau(A) = 0$. \\
~e) The dual abelian variety $A^{\vee}$ is ordinary. \\
~f) The Frobenius map ~$F : H^1(A,{\cal O}_A) \to H^1(A,{\cal O}_A)$ is bijective. \\
 \end{tabular}} \\
 
\vskip 0.1cm

For surfaces the requirement that $Pic^0_X$ be ordinary is clearly of birational nature. As an elementary generalization, one has the following lemma which asserts that ordinarity of the Picard scheme is preserved
under certain group actions. \\

\textbf{Lemma 2 :} \textit{Let $X$ be a surface with ordinary $Pic^0_X$, $Y'= X/G$ be the quotient of $X$ by a
finite group of order prime to $p$ and $Y$ be a desingularization of $Y'$. Then $Pic^0_Y$ is ordinary.} \\

\textbf{Proof :} \\

Let $n = |G|$ be the degree of the quotient map $\pi : X \to Y'$. Since $n$ is prime to $p$,
$\frac{1}{n}\mbox{Trace} : \pi_*{\cal O}_X \to {\cal O}_{Y'}$ defines a splitting of $\pi_*{\cal O}_X$ and we obtain
$$H^1(Y,{\cal O}_Y) = H^1(Y',{\cal O}_{Y'}) \subset H^1(X, {\cal O}_X).$$
Therefore, since all Bockstein operations $\beta_r : H^1(X, {\cal O}_X) \to H^2(X, {\cal O}_X)$ vanish by assumption,
$\beta_r : H^1(Y, {\cal O}_Y) \to H^2(Y, {\cal O}_Y)$ are all zero. Thus we see that $Pic^0_Y$ is smooth. \\

The albanese map $Y \to Alb_Y$ factors through $Y'$. Hence we obtain a surjective homomorphism 
$Alb_X \to Alb_Y$. Therefore, since $Alb_X$ is ordinary so is $Alb_Y$ and the result follows from
Lemma 2(e). $\Box$ \\

\vskip 0.2cm

\textbf{Remarks :}

\begin{itemize}
	
	\item If $p | |G|$, then $Pic^0_Y$ need not be reduced (classical example; Igusa surfaces ([7])). However, the conclusion of Lemma 2 holds for $(Pic^0_Y)_{red}$ of a quotient $Y = X/G$ where $X$ is as in Lemma 2, $G$ is a local commutative group scheme (respectively
a reduced finite group scheme)
with reduced (respectively local commutative) dual $D(G)$ and $G$ acts on $X$ freely. This follows from the exact sequence
$$0 \to D(G) \to Pic_Y \to (Pic_X)^G \to 0$$
where $(Pic_X)^G$ is the fixed point scheme of the action of $G$ on $Pic_X$
(\cite{Je}, Thm.2.1).
	
	\item In (\cite{R1}) Raynaud shows that for any curve there exists an etale Galois cover with galois group of order prime to characteristic, which is non-ordinary. Using this result, one can construct non-ordinary étale Galois covers of (at least) isotrivial fibrations (\cite{JJ2}). Thus we see that the converse of Lemma 2 does not hold.

\end{itemize}

\vskip 0.2cm
It is well known that
a curve $Y$ is ordinary if and only if the jacobian $J_Y$ is ordinary. For varieties of dimension $n \ge 2$, in general ordinarity of $X$ does not imply smoothness of $Pic_X$; for instance, one has singular Enriques surfaces $X$ in characteristic two, which are ordinary with $Pic^0_X = \mu_2$. Therefore, it is more
natural to investigate the relation between ordinarity of $X$ and of the reduced scheme $(Pic_X^0)_{red}$
(or equivalently, the ordinarity of $\mbox{Alb}_X$ (Lemma 2 (e))). \\

\textbf{Lemma 3 :} \emph{ For a smooth projective variety $X$ we have : \\
(i) If $X$ is ordinary in dimension one, then  $(Pic^0_X)_{red}$ is ordinary. \\
(ii) If $Pic^0_X$ is ordinary, then $X$ is ordinary in dimension one.} \\

\textbf{Proof :} \\

The result follows from the identifications
$$H^0(X, B^1) = \{w \in H^0(X, \Omega_X) : dw = 0 = Cw\} = H^1(X_{fl}, \alpha_p)$$
$$\mbox{and} ~H^1(X_{fl}, \alpha_p) = Hom(\alpha_p, Pic^0_X). ~\Box$$ \\

\vskip 0.1cm

\textbf{Lemma 4 : } \textit{Let $S$ be a scheme over $k$ and ${\cal X} \to S$ be a smooth projective map. Then the set
$U = \{s \in S : Pic^0_{X_s} ~\mbox{is ordinary}\}$ is open.} \\

\textbf{Proof :} \\

Let $S' \subset S$ be the open subscheme over which $Pic_{x/S}$ is smooth and apply
(\cite{I2}, Prop. 1.2(b)) to the open subscheme $Pic^0_{\cal X'/S'}$ of the Picard scheme
$Pic_{{\cal X}/S}$.  $\Box$ \\

\vskip 0.1cm

We note that requiring $U = S$ is a stringent condition (cf. \cite{MB}, Chapitre XI, Theoréme 5.2). For instance, an ordinary family of jacobians $J_{X/S} \to S$ of relative dimension $d \ge 2$ over a
smooth projective base curve is isotrivial; this follows from the fact that the corresponding family of ordinary curves $X \to S$ is isotrivial (\cite{S2}, Thm. 5). \\

\vskip 0.2cm

\textbf{II. The main result} \\

\textbf{Theorem 1 :} \textit{Let $\pi : X \to C$ be a smooth isotrivial fibration,
$F$ be the fiber and $G$ be the monodromy group. Suppose that} \\

\begin{tabular}{l}
(i) $F \to F/G$ \textit{is tamely ramified or} \\
(ii) $G$ \textit{is cyclic and the inertia groups are of order} $p^ne, ~n \le 2, ~(p,e) = 1$. \\
\end{tabular} \\

\textit{Then $\pi$ lifts to a fibration ${\cal X} \to {\cal C}$ over $W(k)$ (resp. $W(k)[\zeta]$ where 
$\zeta ^{p^2} = 1$) if (i) (resp. (ii)) holds.} \\

\textbf{Proof :} \\

Let $C' \to C$ be the Galois $G$-cover such that 
$$X \times_C C' = F \times C'.$$

With the restrictions as given in the statement of the lemma, the curves $F, C'$ lift
to curves ${\cal F}, {\cal C}'$  with $G$-action (\cite{GM} for (ii)). Hence  the fibration
$${\cal X} = ({\cal C}' \times {\cal F})/G \to {\cal C}'/G = {\cal C}$$
gives a lifting of $\pi$. $\Box$

\vskip 0.5cm

We apply Theorem 1 to smooth ordinary fibrations (necessarily  isotrivial (\cite{S2}, Thm. 5)) using
the results on the relation between the $p$-rank of a smooth projective curve $F$ (that is the $p$-rank of the jacobian $J_F$)
and the structure of the group of automorphisms $Aut(F)$ of $F$ (\cite{N}), to obtain the following result. \\

\textbf{Corollary 1 :}\textit{ Let $\pi : X \to C$ be a smooth ordinary fibration with fiber $F$ and monodromy group $G$.
Suppose that either $g(F) \le p-2$ or $G$ is cyclic. Then $\pi$ lifts to characteristic zero.} \\

\textbf{Proof :} 

If $g(F) \le p-2$, then $|Aut(F)|$ is not divisible by $p$ (\cite{N}, Corollary on p.598) and Theorem 1 (i) applies. \\

Next suppose that $G$ is cyclic. Since $F$ is ordinary, if $\sigma \in Aut(F) $ is of order $p^n$ and leaves one point of $F$ fixed, it is of order
$p$ (\cite{N}, Corollary on p.599). Therefore, Theorem 1 (ii) applies. $\Box$

\vskip 3cm

\begin{tabular}{ll}

Department of Mathematics & ~~~~~Department of Mathematics \\
Karatekin University & ~~~~~Middle East Technical University \\
Çankırı, Turkey & ~~~~~06800 Ankara, Turkey \\
& \\
arithhesab@gmail.com & \\

\end{tabular}

\end{document}